\documentclass[11pt]{article}
\def\br{\hbox{\it I\hskip -2pt R}}
\def\bn{\hbox{\it I\hskip -2pt N}}
\def\bz{\hbox{\it Z\hskip -4pt Z}}
\def\bq{\hbox{\it l\hskip -5.5pt Q}}

\def\demo{\noindent{\bf Proof.}}
\newtheorem{theorem}{Theorem}[section]
\newtheorem{lemma}[theorem]{Lemma}
\newtheorem{remark}[theorem]{Remark}
\newtheorem{Corollary}[theorem]{Corollary}
\newtheorem{definition}[theorem]{Definition}
\newtheorem{proposition}[theorem]{Proposition}
\begin{document}
\begin{center}
\uppercase{{\bf Stanley-Reisner rings and the radicals of lattice ideals}}
\end{center}
\advance\baselineskip-3pt
\vspace{2\baselineskip}
\begin{center}
{\sc Anargyros Katsabekis, Marcel Morales, and Apostolos Thoma}
\end{center}
\vskip.5truecm\noindent
\vskip.5truecm\noindent
\begin{abstract}
In this article we associate to every lattice ideal $I_{L,\rho }\subset K[x_1,\dots ,x_m]$ a cone
$\sigma $ and a graph $G_{\sigma }$ with vertices the minimal
generators of the Stanley-Reisner ideal of  $\sigma $. To
every polynomial $F$ we assign a subgraph $G_{\sigma }(F)$ of the graph $G_{\sigma }$.
 Every expression of the radical of $I_{L,\rho }$, as a radical of an ideal 
generated by some polynomials $F_1,\dots ,F_s$
 gives a spanning subgraph of $G_{\sigma }$, the $\cup_{i=1}^s G_{\sigma }(F_i)$. 
This result provides a lower bound for the minimal number of generators of $I_{L,\rho }$
 and 
therefore improves the generalized Krull's principal
ideal theorem for lattice ideals. But mainly it provides lower bounds  for the
binomial arithmetical rank and the $A$-homogeneous arithmetical rank of a lattice ideal.
Finally we show, by a family of examples, that the bounds given are sharp.
\end{abstract}

\section{Introduction}

Lattice ideals arise naturally in problems from diverse areas of mathematics, including
toric geometry, integer programming,  dynamical systems, computer algebra, 
graph theory, hypergeometric differential
equations, mirror symmetry and computational statistics, see \cite{E-S}, 
\cite{SST}, \cite{St}, \cite{V}. A fundamental problem in the theory of 
lattice ideals is to determine minimal generators
of the lattice ideal $I_L$ or of the lattice ideal $I_L$ up to radical. 
The main Theorem of this article provides 
a lower bound for the minimal number 
of generators of a lattice ideal,  but also it provides lower bounds  for the
binomial arithmetical rank and the $A$-homogeneous arithmetical rank of a lattice ideal.
The lower bounds depend only on the geometry of the cone associated to the lattice ideal.\\
\par Let $K$ be an algebraically closed field of characteristic  $p\geq 0$. 
A {\em lattice} is a finitely generated free abelian group.
 A {\em partial character} $(L,\rho )$ on $\bz^m$ is a 
homomorphism $\rho $ from a sublattice $L$ of 
$\bz^m$ to the multiplicative group $K^*=K-\{ 0 \}$.
 Given a partial character  $(L, \rho )$ on $\bz^m$, we define 
the ideal  
$$I_{L,\rho }:=(\{{\bf x}^{{\bf \alpha }_+}-
\rho({\bf \alpha } ){\bf x}^{{\bf \alpha } _-}| {\bf \alpha } 
={\bf \alpha } _+-{\bf \alpha} _- \in L\})\subset K[x_1,\dots ,x_m]$$
called {\em lattice ideal}.  Where ${\bf \alpha} _+\in \bn^m$ and  
${\bf \alpha } _-\in \bn^m$ denote the 
positive and negative part of ${\bf \alpha } $,
respectively, and ${\bf x}^{{\bf \beta }}=x_1^{b_1}\cdots x_m^{b_m}$ 
for ${\bf \beta } =(b_1,\dots ,b_m)\in \bn^m$. Lattice ideals are binomial ideals. The 
theory of binomial ideals were developed by Eisenbud and Sturmfels in \cite{E-S}.\\
If $L$ is a sublattice of $\bz^m$, then the {\em saturation} of $L$ is the lattice 
$$ Sat(L):=\{  {\bf \alpha } \in \bz^m | d {\bf \alpha } \in L \ 
for \ some \ d \in \bz^*\}.$$
We say that the lattice $L$ is {\em saturated} if $L=Sat(L)$. 
The lattice ideal $I_{L,\rho }$ is prime if and only if $L$ is saturated.
A prime lattice ideal is called a {\em toric ideal}, while the set of zeroes in $K^m$ 
is an {\em affine toric variety} in the sence of \cite{St}, since we do not require normality. \\

\par Throughout this paper we assume that $L$ is a non-zero  {\em positive sublattice}
 of $\bz^m$, that is $L\cap \bn^m=\{{\bf 0}\}$. This means that the 
lattice ideal $I_{L,\rho }$ is
homogeneous with respect to some positive grading.\\
The  group $\bz^m/Sat(L)$ is free abelian, therefore is isomorphic to $ \bz^n$,
 where $n=m-rank(L)$. 
Let $\psi $ be the above isomorphism, ${\bf e}_1,\dots ,{\bf e}_m$ the unit 
vectors of $\bz^m$  and  
$\psi ({\bf e}_i+Sat(L))={\bf a}_i\in \bz^n$ for $1\leq i\leq m$. \\
Let $A=\{{\bf a}_i| 1\leq i\leq m \}$, we associate to the lattice ideal $I_{L,\rho }$ 
 the rational polyhedral cone  
$\sigma =pos_{\bq }(A):=\{l_1{\bf a}_1+\cdots +l_m{\bf a}_m| l_i\in \bq \ and \ l_i\geq 0\}$.  A cone $\sigma $ is {\em strongly convex} 
if $\sigma \cap -\sigma =\{ {\bf 0} \}$. 
The condition that the lattice $L$ is positive, is equivalent with the condition 
that the cone $\sigma $
is strongly convex.\\

\par We grade $K[x_1,\dots ,x_m]$ by setting $deg_A(x_i)={\bf a}_i$ for $i=1,\dots ,m$.
 We define the {\em $A$-degree } of
the monomial ${\bf x^u}$  to be $$deg_A ({\bf x^u}) :=
u_1{\bf a}_1+\cdots +u_{m}{\bf a}_m\in {\bn} A,$$ where ${\bn} A $ is the semigroup
generated by $A$. 
The lattice ideal $I_{L, \rho }$ is $A$-homogeneous as well as all lattice ideals 
with the same saturation. 
 The {\em
binomial arithmetical rank}
of a binomial ideal $I$ (written $bar(I)$) is the smallest integer $s$
for
which there exist binomials $f_1,\dots ,f_s$ in $I$ such that
$rad(I) = rad(f_1,\ldots ,f_s)$. Hence the binomial arithmetical
rank is an upper bound for the {\em arithmetical rank} of a binomial
ideal (written $ara(I)$), which is the smallest integer $s$
for
which there exists $f_1,\ldots ,f_s$ in $I$ such that
$rad(I) = rad(f_1,\ldots ,f_s)$. Especially, when $I$ is $A$-homogeneous and all the 
polynomials $f_1,\ldots,f_s$ are $A$-homogeneous, the smallest integer 
$s$ is called {\em $A$-homogeneous arithmetical rank} of $I$, 
denoted by $ara_A(I)$. From the definitions, the generalized Krull's
principal ideal theorem and the graded version
of Nakayama's Lemma we deduce the following inequality for a
lattice 
ideal $I_{L,\rho }$:
$$h(I_{L, \rho})\leq ara(I_{L, \rho})\leq ara_A(I_{L, \rho}) 
\leq bar(I_{L, \rho})\leq \mu(I_{L, \rho}).$$
Here $h(I)$
denotes the height and $\mu(I)$ denotes the minimal number of
generators of an ideal $I$. When $h(I) = ara(I)$ the ideal $I$ is called a
{\em
set-theoretic complete intersection} and when $h(I) = \mu(I)$ it
is called a {\em complete intersection}. In several cases the 
lower bound $h(I_{L, \rho})$ given by the generalized Krull's
principal ideal theorem can be improved by using local cohomological methods, 
see \cite{B-S}, \cite{H1}. \\

\par The computation of the numbers $ara(I_{L, \rho}), ara_A(I_{L, \rho}),
 bar(I_{L, \rho})$ is usually an extremely difficult problem
and remains open even for some very simple lattice ideals, like the ideal
of the Macaulay curve $(t^4,t^3u,tu^3,u^4)$ in the three
dimensional projective space, see \cite{E}. 
In the case that we can compute good generating sets for the ideal, 
sharp lower bounds for these numbers may help us to determine the exact 
value of them,  see section 5. 
The numbers $ara(I_{L, \rho})$, $bar(I_{L, \rho})$ and 
$ara_A(I_{L, \rho})$, in the cases that were known up to this work, 
were either identical 
or very close to each other, see for
example
\cite{B-M-T1}, \cite{B-M-T2}, \cite{Eto}, \cite{R-V}, \cite{Th}. 
Also, there was no known example of a lattice ideal $I_{L, \rho}$ with the property
$ara(I_{L, \rho})\not=ara_A(I_{L, \rho})$. 
In this work, by providing good lower bounds for $ara_A(I_{L, \rho})$ and  $bar(I_{L, \rho})$ 
and using the result of Eisenbud, Evans and Storch, see \cite{E-E} and \cite{Sto},
that $ara(I_{L, \rho})$ is bounded above by  the dimension $m$ of the space $K^m$,
we show that there can be 
very large differences beetween these numbers. 
For example, using the results of section 5  and putting  
$n=10$ we have an example of a lattice ideal for which  the height is equal to
80, the $ara(I_{L, \rho})$ is smaller than 90 by Eisenbud, Evans and Storch, 
while the $ara_A(I_{L, \rho})$ is exactly 1740 and $bar(I_{L, \rho})$ is exactly 1860.\\  

\par In section 2 we recall some basic facts about lattice ideals, which
are necessary for the formulation and proof  of the main Theorem 4.1.\\
In section 3 we introduce a graph $G_{\sigma }$ with vertices the minimal generators 
of the Stanley-Reisner ideal  of the cone $\sigma $ associated to the lattice ideal. \\
In section 4 we state and prove the main theorem of the article, Theorem 4.1, which
 provides lower bounds 
for the A-homogeneous arithmetical rank, the binomial arithmetical rank and the
minimal number of generators of a lattice ideal.\\
In section 5 we compute these  bounds for a special
class of lattice ideals. In this case, we show that the lower bounds given by Theorem
4.5 cannot be improved, by computing the exact value of the 
A-homogeneous arithmetical rank and the binomial 
arithmetical rank for certain lattice ideals.

\section{Basics on Lattice ideals }

Let $L$ be a nonzero positive sublattice of $\bz^m $ and $(L, \rho )$ be a partial
character on $\bz^m $.
\begin{definition} If $p$ is a prime number, we define $Sat_p(L)$ and $Sat'_p(L)$ 
to be the largest sublattices of $Sat(L)$
 containing $L$ such that $Sat_p(L)/L$ has order a power of $p$ and 
$Sat'_p(L)/L$ has order relatively prime to $p$.
If $p=0$, we define $Sat_p(L)=L$ and $Sat'_p(L)=Sat(L)$.
\end{definition} 
\begin{theorem}\cite{E-S}  Let $(L, \rho )$ be a partial
character on $\bz^m $. Write $g$ for the order of $Sat'_p(L)/L$. There are $g$ distinct
characters $\rho _1,\dots ,\rho _g$ of $Sat'_p(L)$ extending $\rho $ and for each $j$ 
a unique character $\rho '_j$ of $Sat(L)$ extending $\rho _j$. There is a unique 
partial character $\rho '$ of $Sat_p(L)$ extending $\rho $. The radical,
associated primes and minimal primary decomposition of 
$I_{L,\rho }$ are:\\
$$rad(I_{L,\rho })=I_{Sat_p(L), \rho '},$$
$$ Ass(K[x_1,\dots ,x_m]/I_{L,\rho })=\{I_{Sat(L),\rho '_j}| j=1,\dots ,g\}$$
and $$I_{L,\rho }=\bigcap^g_{j=1}I_{Sat '_p(L),\rho _j}$$
where $I_{Sat '_p(L),\rho _j}$ is $I_{Sat(L),\rho '_j}$-primary. 
In particular, if $p=0$, then $I_{L, \rho }$ is 
a radical ideal. The associated primes $I_{Sat(L),\rho '_j}$ 
of $I_{L, \rho }$ are all minimal
and have the same codimension $rank(L)$.
\end{theorem}
\par We decompose the affine space $K^m$ into $2^m$ {\em coordinate
cells}, $$(K^*)^E:=\{(q_1, \ldots ,q_m)\in K^m | q_i \not= 0 \ for
\ i\in E, q_i=0\ for \ i\notin E \},$$
where $E$ runs over all subsets of $\{ 1,\ldots ,m \}$. 
 We denote by $K[E]:=K[\{x_i| i\in E\}]$. 
Let $P=(x_1,\dots ,x_m)$ be a point of $K^m$ then  
$$P_E:=(\delta^{E}_1x_1, 
\delta^{E}_2x_2, \ldots ,
\delta^{E}_mx_m)\in K^m, $$ where $\delta^{E}_i=1$
if $i\in E$ and $\delta^{E}_i=0$
if $i\notin E$. Note that if $P\in (K^*)^{ \{ 1,\ldots ,m \} }$
then $P_E\in (K^*)^E$. \\ 
A face $F$ of ${\sigma}$ is any set of the form $$F= \sigma \cap
\{{\bf x}\in \br ^n: {\bf cx} = 0\}, $$
where ${\bf c}\in \br^n$ and ${\bf cx} \ge 0$ for all points ${\bf x}\in \sigma$.
Faces of dimension one are called {\em extreme rays}.  
If the number of the extreme rays of a cone coincides with
the dimension (i.e. the extreme rays are linearly independent), the
cone is called {\em simplex cone}.\\
Let $S$ be a subset of the cone $\sigma $, then 
 $E_S:=\{i \in \{1, \dots ,m\}|{\bf a}_i\in  S\}$.
To simplify the notation we denote the point $P_{E_S}$ by $P_S$ and
the cell $(K^*)^{E_S}$ by $(K^*)^S$.
The $n$-dimensional algebraic torus $( K^*)^n$ acts
on the affine
$m$-space $ K^m$ via $$(x_1, \ldots ,x_m) \rightarrow (x_1{\bf t}^{{\bf a}_1},
\ldots ,x_m{\bf t}^{{\bf a}_m}).$$
\par Let $j\in \{1,\dots ,g\}$. 
The affine toric variety ${\bf X}_{A,j}:=V( I_{Sat(L),\rho '_j}) $ 
is the Zariski-closure of the
$( K^*)^n$-orbit of a point $P_{j}=(c_{j1},c_{j2},\ldots ,c_{jm})$, where all $c_{ji}$ 
are different from zero. Note that
the ideal  $I_{Sat(L),\rho '_j}$ is the kernel of the $K$-algebra homomorphism
$$\phi_j: K[x_1,\ldots, x_m]\rightarrow
 K[t_1,\ldots,t_{n},t_1^{-1},\ldots ,t_{ n }^{-1}]$$ 
 given by $$\phi_j(x_i) = c_{ji}{\bf t}^{{\bf a}_i }
 \qquad \mbox{for
 all }i = 1,\ldots,m. $$ 
The $( K^*)^n$-orbits
on the  affine toric variety ${\bf X}_{A,j}$ are in order-preserving
bijection with the faces of the cone $\sigma $, 
see \cite{F}, \cite{I}, \cite{O}, for every
$j$.
Note that our cone $\sigma $ is the dual of the cone that is used
to define the toric variety in the above references.\\ 
 Actually the orbit
corresponding to the face $F$ is the 
orbit of the point $(P_j)_F$ and the toric variety is the disjoint 
union of the orbits of the points 
$(P_j)_F$, for all the faces $F\in \sigma $, i.e.
$${\bf X}_{A,j}=\cup_{F\in \sigma }O((P_j)_F).$$ Each orbit $O((P_j)_F)$ corresponds to the relative interior of the face $F$. The orbit
$O((P_j)_F)$ is in the cell $(K^*)^F$
and there are no points of the toric varieties ${\bf X}_{A,j}$ that are in the cells 
$(K^*)^{E}$, where $E$ is not in the form $E _F$ for 
a face $F$ of $\sigma $.\\
From the theorem 2.2 we have $V(I_{L,\rho })=\cup _{j=1}^{g}{\bf X}_{A,j}$. Therefore $V(I_{L,\rho })$ has points only on the cells in the form $(K^*)^F$
for some face $F$ of the cone $\sigma $.

\section{ Stanley-Reisner rings }

Given a set $Y\subset \bz^n$,  the set of all nonnegative linear combinations
$x=\l_1{\bf y}_1+\cdots +l_s{\bf y}_s$, where ${\bf y}_1,\ldots ,{\bf y}_s\in Y$, 
$l_1,\ldots ,l_s\in \bq$, is called the {\em positive hull} of $Y$, $pos_{\bq}(Y)$.\\
\par Let $\sigma \subset \bq ^n$ be a strongly convex rational polyhedral cone 
and let ${\bf R}_{\sigma }=\{{\bf r}_1, \ldots , {\bf r}_t\}$ a set 
of integer vectors, one for each extreme ray of $\sigma $, therefore 
$\sigma =pos_{\bq}({\bf r}_1,\ldots , {\bf r}_t )$. The vectors ${\bf r}_i$
are called {\em extreme vectors} of $\sigma $. We consider 
the polynomial ring $K[Y_1,\ldots,Y_t]$ by
taking one variable $Y_i$ for each vector ${\bf r}_i$.
Let $M=Y_{i_1}^{n_1}\cdots Y_{i_l}^{n_l}$ be a monomial, we shall 
denote by $pos_{\bq}(M)$ the positive hull of the vectors 
${\bf r}_{i_1},\ldots , {\bf r}_{i_l}$. \\
The {\em relative interior} of $M$, 
$relint(M)=relint_{\bq}({\bf r}_{i_1},\ldots ,{\bf r}_{i_l}),$ 
is the set of all positive rational linear combinations 
of ${\bf r}_{i_1},\ldots ,{\bf r}_{i_l}$.\\
The {\em boundary} of $M$ is defined to be 
$$\partial (M)= \partial ({\bf r}_{i_1},\ldots ,{\bf r}_{i_l}):=
pos_{\bq}({\bf r}_{i_1},\ldots ,{\bf r}_{i_l})-relint_{\bq}({\bf r}_{i_1},\ldots ,{\bf r}_{i_l}),$$
which is the union of all proper faces of the cone 
$pos_{\bq}({\bf r}_{i_1},\ldots ,{\bf r}_{i_l})$.\\
By $F(M)$ we denote the minimal face of $\sigma $ that 
contains $\{{\bf r}_{i_1},\ldots ,{\bf r}_{i_l}\}$, i.e. 
$$F(M)=\cap_{\{{\bf r}_{i_1},\ldots ,{\bf r}_{i_l}\}\subset F}F,$$ 
since any intersection of faces of $\sigma $  is a face of $\sigma $.\\
The  {\em Stanley-Reisner ring} of $\sigma $ is the K-algebra 
$$K[\sigma ]=K[Y_1,\ldots ,Y_t]/I_{\sigma },$$
where $I_{\sigma } $ is the  {\em Stanley-Reisner ideal} generated by all 
square-free monomials $M=Y_{i_1}Y_{i_2}\cdots Y_{i_l}$ such that 
$pos_{\bq}(M)$ is not a face of $\sigma $. \\
\par The ideal $I_{\sigma } $ is a monomial ideal, so there is a unique  set 
$\{M_1,\ldots ,M_q \}$ of minimal
square-free monomial generators of $I_{\sigma } $. \\
\begin{definition} We associate to the cone $\sigma $ a graph $G_{\sigma }$ with vertices the set 
$\{M_1,\dots ,M_q\}$of minimal
monomial generators of $I_{\sigma} $. There is an edge between the vertices
$M_i$ and $M_j$ if $relint_{\bq}(M_i)\cap relint_{\bq}(M_j)\not= \{ 0 \}$.\\
\end{definition}
\begin {remark}{\rm $G_{\sigma }=\emptyset$ if and only if $\sigma $ is simplex cone.}  
\end {remark}

The next Theorem gives an equivalent condition for a square-free monomial to be
minimal generator of $I_{\sigma}$.\\
\begin{theorem} The monomial $M$ is a minimal generator of $I_{\sigma}$ iff \\
i) for every proper divisor $N$ of $M$,  $pos_{\bq}(N)$ is a face of $\sigma $\\
ii) the positive hull $pos_{\bq}(M)$ is a proper subset of $F(M)$ \\
iii) $pos_{\bq}(M)$ is a {\em simplex cone} and every proper face of $pos_{\bq}(M)$
is a face of $\sigma $. \end{theorem}
\demo Suppose that $M=Y_{i_1}Y_{i_2}\cdots Y_{i_l}$ 
is a minimal generator of $I_{\sigma}$.\\
i) Assuming that $pos_{\bq}(N)$ is not a face of $\sigma $ 
we have $N \in I_{\sigma }$ from the definition
of Stanley-Reisner ideal. But this contradicts the fact 
that $M$ is a minimal generator of $I_{\sigma}$.\\
ii) The positive hull of $M$ is not a face of $\sigma$, while $F(M)$ 
is a face of $\sigma$. Thus $pos_{\bq}(M)\not=F(M)$ and certainly $pos_{\bq}(M)\subset F(M)$.\\
iii) Assume that ${\bf r}_{i_1},{\bf r}_{i_2},\ldots ,{\bf r}_{i_l}$ 
are not linearly independent and consider a linear relation 
$d_{i_1}{\bf r}_{i_1}+d_{i_2}{\bf r}_{i_2}+\cdots +d_{i_l}{\bf r}_{i_l}=0$ 
between them, with atleast one $d_{i_j}\not=0$. Then, 
since $\sigma$ is strongly convex, there will be positive and negative
 coefficients $d_{i_j}$ in the previous relation.  Let $P$ be the  subset of $\{i_1,\ldots,i_l\}$ 
consisting from all indices $i_j$, such that the corresponding $d_{i_j}$ 
is positive. Then $P$ is not empty and
proper. 
Therefore $N=\Pi _{i \in P}Y_i$ is a proper divisor of $M$ 
which means that $pos_{\bq}(N)$ is a face $F$ of $\sigma $. 
Let ${\bf c}_F$ be a vector defining the face $F$. Considering the
dot product of ${\bf c}_F$ and 
$d_{i_1}{\bf r}_{i_1}+d_{i_2}{\bf r}_{i_2}+\cdots +d_{i_l}{\bf r}_{i_l}$ 
we have a contradiction, namely a negative number equal to zero. Therefore
$pos_{\bq}(M)$ is a simplex cone.\\
Let $F$ be a proper face of $pos_{\bq}(M)=pos_{\bq}({\bf r}_{i_1},\ldots , {\bf r}_{i_l})$.
Then $F=pos_{\bq}({\bf r}_{j_1},\ldots , {\bf r}_{j_k})$, where  
$\{{\bf r}_{j_1},\ldots , {\bf r}_{j_k}\}$ is aproper subset of $\{{\bf r}_{i_1},\ldots , {\bf r}_{i_l}\}$. Then $N=Y_{j_1}\cdots Y_{j_q}$ is a
proper divisor of $M$, therefore $F$ is a face of $\sigma $.\\
Suppose that i), ii) and iii) are true. Then ii) give us that $M$ is a generator of the 
Stanley-Reisner ideal, while i) ensure that $M$ is minimal.
\bigskip
\newline
\par The following lemma will be usefull in the proof of the theorem 4.1.

\begin{lemma}The monomial ${\bf x^u}\in K[E_F] $ iff $deg_A({\bf x^u})
\in F$.
\end{lemma}
\demo Obviously, ${\bf x^u}$ belongs to $K[E_F] $ 
implies that $deg_A({\bf x^u})$ is in $F$. Suppose that 
$$deg_A({\bf x^u})= u_1{\bf a}_1+\cdots +u_{m}{\bf a}_m \in F.$$ 
Then $$0={\bf c}_F (\sum _{i=1}^{m} u_i{\bf a}_i)=\sum _{i=1}^{m}u_i{\bf c}_F {\bf a}_i,$$ 
where ${\bf c}_F$ is any vector that
defines the face $F$. All the terms ${\bf c}_F
{\bf a}_i$ are non-negative and every $u_i\geq 0$, therefore we have that 
$u_i=0$ whenever ${\bf c}_F
{\bf a}_i$ is positive. Thus ${\bf x^u}\in K[E_F] $.

\section{Radical of a Lattice ideal}

We consider a lattice ideal $I_{L,\rho }\subset K[x_1,\ldots ,x_m]$ and 
the strongly convex rational polyhedral cone
$\sigma =pos_{\bq}(A)\subset  {\bq }^n$  corresponding to $I_{L,\rho }$. Let 
$I_{\sigma }\subset K[Y_1,\dots ,Y_t]$ be the
Stanley-Reisner ideal of the cone $\sigma $, where $t$ is the number of extreme rays
of the cone $\sigma $.
Let $N=x_{i_1}^{n_1}\cdots x_{i_s}^{n_s}$ be a monomial in $K[x_1,\ldots ,x_m]$. 
Set $A_N:=\{{\bf a}_{i_1},\ldots, {\bf a}_{i_s}\}$, we define the cone of $N$ 
to be $$cone(N):=\cap _{A_N\subset pos_{\bq}({\bf r}_{j_1},\ldots ,{\bf r}_{j_l})}
pos_{\bq}({\bf r}_{j_1},\ldots ,{\bf r}_{j_l})\subset \sigma.$$ 
 Note that $pos_{\bq}(A_N)\subset cone(N)$. Also, the  $cone(N)$ is not necessarily in the form $$pos_{\bq}({\bf r}_{j_1},\ldots ,{\bf r}_{j_t})$$
for some extreme vectors ${\bf r}_{j_1},\ldots ,{\bf r}_{j_t}$ of $\sigma $.
 But in the case that 
every one of ${\bf a}_{i_1},\ldots, {\bf a}_{i_s}$ 
belongs to some  extreme ray of $\sigma $, we have that $cone(N)=
pos_{\bq}({\bf a}_{i_1},\ldots, {\bf a}_{i_s})$. \\
\par Let $F$ be a polynomial in $K[x_1,\ldots ,x_m]$, we associate to $F$ the  
induced subgraph $G_{\sigma }(F)$ of $G_{\sigma }$ with vertices those 
$M_i$ with the property that there 
exist a monomial $N$ in $F$ such that $cone(N)=pos_{\bq}(M_i)$. The 
{\em induced subgraph} of a graph $G$
by certain vertices $V$ is the subgraph of $G$ with these vertices and edges those edges of $G$
that have both vertices in $V$.\\
A subgraph $H$ of a graph $G$ is called a {\em spanning subgraph} if $V(H)=V(G)$, where $V(G)$ denotes the set of 
vertices of a graph $G$.\\

\begin{theorem} Every expression of $rad(I_{L,\rho })=rad(F_1,\ldots ,F_s)$  gives a
spanning subgraph of $G_{\sigma }$, the $\cup_{i=1}^s G_{\sigma }(F_i)$. 
\end{theorem}
\demo  Suppose that $rad(I_{L,\rho })=rad(F_1,\ldots ,F_s)$ and 
let $M=Y_{i_1}Y_{i_2}\cdots Y_{i_l}\not= 0$ be a minimal generator of the 
Stanley-Reisner ideal of $\sigma $. We will prove that there exists a monomial 
$N=x_{i_1}^{n_1}\cdots x_{i_s}^{n_s}$ in some $F_i$ such that $cone(N)=pos_{\bq}(M)$.\\
Let us consider the point $(P_j)_{\partial (M)}$, for any $j\in \{1,\dots ,g\}$. We
divide the proof into three steps:\\
{\bf $1^{st}$ step}. 
We claim that $(P_j)_{\partial (M)}$ is not a point of $V(I_{L,\rho })$.
Recall that $(P_j)_{\partial (M)}$ belongs to the cell $(K^*)^{\partial(M)}$.
But since every point of $V(I_{L, \rho })$ belongs to a cell $(K^*)^F$ for some face $F$ of
$\sigma $, 
it is enough to prove that $E _{\partial(M)}$ is not 
in the form $E_F$ for 
a face $F$ of $\sigma $. Suppose that $E_{\partial(M)}=E_F$
for a face $F$ of $\sigma $. \\
Note that $M\not=0$ is a minimal generator of the Stanley-Reisner ideal of the cone $\sigma $ and therefore $dim(pos_{\bq}(M))\ge 2$. Also $\sigma =pos_{\bq}(A)$ which implies that
for every  extreme vector ${\bf r}_k$ of $\sigma $ there exist ${i_k}\in \{1,\dots ,m\}$
such that ${\bf a}_{i_k}=\lambda {\bf r}_k$, for some $\lambda \in \bq$.  
Then ${\bf r}_k\in pos_{\bq}(M)$ iff ${\bf r}_k\in {\partial(M)}$ iff 
 ${\bf a}_{i_k}\in {\partial(M)} $ iff ${i_k}\in E_{\partial(M)}$.
Also ${i_k}\in E_F$ iff ${\bf a}_{i_k}\in F$ iff ${\bf r}_k\in F$. 
Therefore $pos_{\bq}(M)=F$, since every face of $\sigma $ is generated by extreme vectors.
But this contradicts the fact that $M$ is a generator of the Stanley-Reisner ideal
and the claim is proved.\\
Therefore
$(P_j)_{\partial (M)}$ cannot be a zero of all the $F_i$. Thus there exists atleast 
one $i$ such that $F_i((P_j)_{\partial (M)})\not=0$.
Let $N$ be a monomial in $F_i$ such that $N((P_j)_{\partial (M)})\not=0$. \\
We have, from the definition of $(P_j)_{\partial (M)}$ and the fact  
$N((P_j)_{\partial (M)})\not=0$, that 
$A_N\subset \partial (M)
\subset pos_{\bq}(M)$. 
Therefore $cone(N)\subset pos_{\bq}(M)$.\\
The last condition implies that $deg _{A}(N)\in pos_{\bq}(M)$, that means either\\ 
$deg _{A}(N)\in relint(pos_{\bq}(M))$ or $deg _{A}(N)\in \partial(M)$. \\
{\bf $2^{nd}$ step}.
We claim that always we can find a monomial $N$ in $F_i$ such that
$deg _{A}(N)\in relint(pos_{\bq}(M))$. \\
Suppose that $deg _{A}(N)\in \partial(M)$. But $M$ 
is a minimal generator of the Stanley-Reisner ideal and
therefore, from theorem 3.3, we have that $deg _{A}(N)$ belongs 
to a face $F$ of the cone $\sigma $ such
that $F\subset \partial(M)$. \\
The polynomial $F_i$ belongs to the lattice ideal $I_{L,\rho }$, which is 
$A$-homogeneous and therefore 
 has a decomposition $F_i=F_{i1}+\cdots +F_{is}$ into $A$-homogeneous components.
By lemma 3.4, $deg _{A}(N)$ belongs to a face $F$ implies that the $A$-component, $F_{ij}$, 
of $N$ belongs to
$K[E_F]$,  since all monomials in $F_{ij}$ have the same
$A$-degree. Note that $F\subset \partial(M)$ therefore 
$((P_j)_{\partial (M)})_F=(P_j)_F$. Thus, since 
$F_{ij}$ involves variables belonging only to the face $F$, we have
$F_{ij}((P_j)_{\partial (M)})=F_{ij}((P_j)_F)=F_{ij}(P_j)$ which is zero
because $P_j \in V(I_{L,\rho })$.\\
But then $F_{ij}((P_j)_{\partial (M)})=0$ and $F_i((P_j)_{\partial (M)})\not= 0$, 
so there exist a different monomial $N'$ in a different
$A$-homogeneous component of $F_i$ such that $N'((P_j)_{\partial (M)})\not=0$. 
This cannot be repeated indefinitely, since
$F_i$ has finitely many $A$-homogeneous components. So we conclude that
there must be an $N$ in $F_i$ such that $deg _{A}(N)\in relint(pos_{\bq}(M))$ and
$N((P_j)_{\partial (M)})\not=0$. \\
{\bf $3^{rd}$ step}. For a set $S\subset \sigma$ we define ${\bf R}_S$ to be the
set of extreme vectors of $\sigma $ that belong to $S$. 
We will show that  a monomial $N$ with the property 
$deg _{A}(N)\in relint(pos_{\bq}(M))$ and
$N((P_j)_{\partial (M)})\not=0$ satisfies $cone(N)=pos_{\bq}(M)$.
Let ${\bf a}_i\in A_N$, then  $N((P_j)_{\partial (M)})\not=0$ implies that 
${\bf a}_i\in {\partial (M)}$. By theorem 3.3 we conclude that ${\bf a}_i\in F$
for some face of $\sigma $. Therefore $F({\bf a}_i)\subset F\subset {\partial (M)}$,
where $F({\bf a}_i)$ denotes 
the smallest face that contains ${\bf a}_i$. We have that
$${\bf R}_{F({\bf a}_i)}\subset {\bf R}_{F}\subset 
{\bf R}_{\partial (M)}={\bf R}_{pos_{\bq}(M)}.$$
Now we claim that if ${\bf a}_i\in pos_{\bq}({\bf R})$, 
for some ${\bf R}\subset {\bf R}_{\sigma }$, then ${\bf R}_{F({\bf a}_i)}\subset {\bf R}$. 
Let ${\bf a}_i=\Sigma_{{\bf r}_i\in {\bf R}}l_i{\bf r}_i$, with  $l_i\ge 0$.
 Multiplying by ${\bf c}_{F({\bf a}_i)}$ 
a vector that defines the face $F({\bf a}_i)$, we have that $l_j=0$ whenever 
${\bf r}_{j}\notin F({\bf a}_i)$. So in fact  
${\bf a}_i=\Sigma_{{\bf r}_i\in {\bf R}_{F({\bf a}_i)}\cap {\bf R}}l_i{\bf r}_i$. Note
also that $pos_{\bq}({\bf R}_{F({\bf a}_i)}\cap {\bf R} )$ is a face of $\sigma $ by
theorem 3.3,
since ${\bf R}_{F({\bf a}_i)}\cap {\bf R}$ is a proper subset of ${\bf R}_{pos_{\bq}(M)}$.
We conclude that 
$$F({\bf a}_i)\subset pos_{\bq}({\bf R}_{F({\bf a}_i)}\cap {\bf R} )
\subset pos_{\bq}({\bf R}_{F({\bf a_i})}) = F({\bf a}_i).$$ 
Therefore ${\bf R}_{F({\bf a}_i)}\cap {\bf R}={\bf R}_{F({\bf a_i})}$
which implies the claim ${\bf R}_{F({\bf a}_i)}\subset {\bf R}$.\\
To prove that
 $cone(N)=pos_{\bq}(M)$  is enough 
to prove that 
$\cup_{{\bf a}_i\in A_N} {\bf R}_{F({\bf a}_i)}=
{\bf R}_{pos_{\bq}(M)}$.
We have just proved that $\cup {\bf R}_{F({\bf a}_i)}\subset {\bf R}_{pos_{\bq}(M)}$. 
If they are not equal
then $\cup {\bf R}_{F({\bf a}_i)}\subset F$, for some face $F$ of $\sigma $, 
since $M$ is a minimal generator
of the Stanley-Reisner ideal. But if $\cup {\bf R}_{F({\bf a}_i)}\subset F$ 
then $deg_{A}(N)\in F$. Which is a contradiction,
since $deg_{A}(N)\in relint(pos_{\bq}(M))$. Therefore we have proved that 
for every minimal generator $M$ of
the Stanley-Reisner ideal  of $\sigma $ there exists atleast one monomial 
$N$ in some $F_i$ such that 
$cone(N)=pos_{\bq}(M)$ and even more, $deg_{A}(N)\in relint(pos_{\bq}(M))$
and $A_N\subset {\partial (M)}$.
 
\begin{theorem} Let $F\in K[x_1,\ldots ,x_m]$ 
 be an $A$-homogeneous polynomial, 
then the graph $G_{\sigma }(F)$ is complete.
\end{theorem}
\demo  Suppose that $G_{\sigma }(F)$ is not empty and that 
$M_i$, $M_j$ are two vertices of $G_{\sigma }(F)$. Let $N_i$ and $N_j$ be the 
corresponding monomials in $F$ with $deg _{A}(N_i)\in
relint_{\bq}(M_i)$ and $deg _{A}(N_j)\in
relint_{\bq}(M_j)$, see the proof of Theorem 4.1. Using the fact that 
$F$ is $A$-homogeneous we get $deg _{A}(N_i)= deg _{A}(N_j)$.  
Thus $relint_{\bq}(M_i)\cap relint_{\bq}(M_j)\not= \{ 0 \}$ and 
therefore, from the definition of $G_{\sigma}$, there is an edge between them.\\
It follows that the subgraph $G_{\sigma }(F)$ is complete, 
since for any two vertices $M_i$ and $M_j$ of  $G_{\sigma }(F)$
 there is an edge between them.
\bigskip
\newline
\par Combining Theorems 4.1 and 4.2 we have the following corollary:

\begin{Corollary} Every expression of $rad(I_{L,\rho })=rad(F_1,\ldots ,F_s)$, 
where each $F_i$ is $A$-homogeneous polynomial, gives a
subgraph of $G_{\sigma }$ which is spanning and is a union of complete 
subgraphs. \\
\end{Corollary}
\medskip
Note that binomials  belonging to $I_{L,\rho }$ are always
A-homogeneous and therefore we have the following corollary: \\

\begin{Corollary} Every expression of $rad(I_{L,\rho })=rad(B_1,\ldots ,B_s)$, 
where each $B_i$ is binomial, gives a
subgraph of $G_{\sigma }$ which is spanning and each binomial 
contributes two vertices and 
 an edge joining them or just one vertex or nothing.  \\
\end{Corollary}
 
For a graph $G$ we denote by $c_G$  the smallest number $s$ of 
complete subgraphs $G_i$ of $G$, such that the subgraph
$\cup_{i=1}^sG_i$ of $G$ is spanning. While by $b_G$ 
we denote the smallest number $s$ of
 subgraphs $B_i$ of $G$,
consisting of two vertices and an edge or just a vertex, such that the subgraph
$\cup_{i=1}^sB_i$ is spanning. 
Then Corollaries 4.3 and 4.4 imply that:
\begin{theorem} For a lattice ideal $I_{L,\rho }$ with associated cone $\sigma $ we have
$c_{G_{\sigma }}\leq ara_{A}(I_{L,\rho })$ and $b_{G_{\sigma }}\leq bar(I_{L,\rho })$.
\end{theorem}
Note that $b_{G_{\sigma }}\ge q/2$, where we recall that $q$
is the minimal number of generators of $I_{\sigma }$, and $c_{G_{\sigma }}$ is greater than
or equal to the number of connected components of the graph $G_{\sigma } $.\\
Also note that the above bounds depend only on the graph $G_{\sigma } $,
 i.e.
lattice ideals with associated cones rationally affine equivalent  have exactly the same bound.
Two cones are called  {\em rationally affine equivalent} if there is a 
rational affine tranformation mapping
the first cone to the second bijectively.
\begin{Corollary} Every expression of $I_{L,\rho }=(B_1,\ldots ,B_s)$, where each $B_i$ is binomial, gives a
subgraph of $G_{\sigma }$ which is spanning and each binomial
 contributes two vertices and 
 an edge joining them or just one vertex or nothing. In particular 
$max\{b_{G_{\sigma }}, h(I_{L,\rho })\} \leq \mu(I_{L,\rho })$. \\
\end{Corollary}
The above Corollary gives a lower bound for the minimal number of generators of $I_{L,\rho }$
which improves the generalized Krull's principal ideal theorem, see also remark 5.6.\\

\section{The lower bounds are sharp}

The aim of this last section is to explicitly compute  the bounds for the $A$-homogeneous
arithmetical rank and the binomial arithmetical rank, obtained from Theorem 4.5,
 for a special class of lattice ideals 
and show that the lower bounds given are sharp. 
This will be done by computing the exact values of the above numbers
 and proving that they are identical with the corresponding bounds, 
for a certain class of lattice ideals.\\
\par We consider the set of vectors 
$A_n=\{2{\bf e}_i+{\bf e}_j:1\leq i,j\leq n,$ $i\not= j\}$,
where $n\ge 2$ and ${\bf e}_i, 1\leq i\leq n,$ is the canonical basis of $K^n$. 
The toric ideal $I_{A_n}$ of $A_n$, see \cite{St}, is the kernel of the 
$K$-algebra homomorphism $\phi : K[\{x_{ij}|1\leq i,j\leq n,i\not= j\}]\rightarrow K[t_1,\ldots,t_{n}]$ given by $$\phi(x_{ij}) = t_i^2t_j.$$
Let $I_{L,\rho }$ be any lattice ideal with associated cone $\sigma =pos_{\bq}(A_n)$ 
or rationally affine
equivalent to the cone $pos_{\bq}(A_n)$. \\
We define the following vectors in ${\bq }^{n}$, with coordinates:
$$({\bf c}_T)_s= \left \{\begin{array} {ll} 0, & for \ s\in T\\
  1, & \mbox{otherwise,} \end{array} \right.$$ 
$$({\bf c}_{i,T})_s= \left \{\begin{array} {lll} -1, & for \ s=i\\
 2, & for \ s\in T\\ 3 & \mbox{otherwise,} \end{array} \right.$$
where $1\leq s \leq n$, $T$ is a subset of $\{1,\dots ,n\}$ 
and $1\leq i \leq n$, $i\notin T$.\\
Note that the $pos_{\bq}(2{\bf e}_i+{\bf e}_j)$ is an extreme ray of the 
cone $pos_{\bq}(A_n)\subset {\bq }^{n}$ with defining vector ${\bf c}_{i,\{j\}}$. Therefore the cone $pos_{\bq}(A_n)$ has $n(n-1)$ extreme rays. 
\par We consider the Stanley-Reisner ideal 
$I_{pos_{\bq}(A_n)}\subset K[\{Y_{ij}|1\leq i,j\leq n,i\not= j\}]$. 
For the graph $G_{pos_{\bq}(A_n)}$ we have the following result: 
\begin{proposition} There are $9\pmatrix{n \cr 3\cr}+12\pmatrix{n \cr 4\cr}$ vertices,
$15\pmatrix{n \cr 3\cr}+18\pmatrix{n \cr 4\cr}$ edges and 
$\pmatrix{n \cr 3\cr}+\pmatrix{n \cr 4\cr}$ connected
 components in the graph  $G_{pos_{\bq}(A_n)}$.
\end{proposition}
\demo  We claim that the minimal generators of   $I_{pos_{\bq}(A_n)}$ 
are  the $3\pmatrix{n \cr 3\cr}$ quadratic monomials in the
form $Y_{ij}Y_{kj}$, the $6\pmatrix{n \cr 3\cr}$ monomials in the form $Y_{ij}Y_{ki}$ 
and the $12\pmatrix{n \cr 4\cr}$ monomials
 in the form $Y_{ij}Y_{kl}$,
 where $i,j,k,l \in \{1,\ldots,n\}$. Here we adopt the convention that $\pmatrix{n \cr k\cr}=0$ for $k>n$. \\ 
The relation $(2{\bf e}_i+{\bf e}_j)+(2{\bf e}_k+{\bf e}_j)=
(2{\bf e}_j+{\bf e}_i)+(2{\bf e}_k+{\bf e}_i)$ 
shows that $pos_{\bq}(Y_{ij}Y_{kj})$ cannot be a face of the cone $pos_{\bq}(A_n)$.
In the contrary case,  taking the dot product with its defining vector in the two parts of the equality 
we get zero equal to a positive number, which is 
a contradiction. Thus 
$Y_{ij}Y_{kj}$ is a
generator of $I_{pos_{\bq}(A_n)}$.
Similarly, the relations 
$2(2{\bf e}_i+{\bf e}_j)+(2{\bf e}_k+{\bf e}_i)=
2(2{\bf e}_i+{\bf e}_k)+(2{\bf e}_j+{\bf e}_i)$ and $(2{\bf e}_i+{\bf e}_j)+(2{\bf e}_k+{\bf e}_l)=
(2{\bf e}_i+{\bf e}_l)+(2{\bf e}_k+{\bf e}_j)$
show that $Y_{ij}Y_{ki}$ and $Y_{ij}Y_{kl}$ are generators of $I_{pos_{\bq}(A_n)}$. 
They are minimal, since there is no linear monomial in $I_{pos_{\bq}(A_n)}$.
Next we show that there is no other minimal generator of the Stanley-Reisner ideal.
The only square free monomials of degree greater 
than or equal to two that are not divided by 
the previous quadratic minimal generators 
are in the form $M_{i,T}=\prod _{j\in T}Y_{ij}$ for some 
$T\subset \{1,\ldots,i-1,i+1,\ldots,n\}$ or $M_{\{i,j\}}=Y_{ij}Y_{ji}$.
But $pos_{\bq}(M_{i,T})$ and $pos_{\bq}(M_{\{i,j\}})$ are faces whose defining
vectors are  ${\bf c}_{i,T}$  and ${\bf c}_{\{i,j\}}$.\\
We define the {\em index} of a $Y_{ij}$ to be the set $\{i,j\}$ and the {\em index of
a monomial} $M\in K[\{Y_{ij}|1\leq i,j\leq n,i\not= j\}]$ to be the union of the indices of the variables in $M$.\\
 Let $M$ and $N$ be minimal generators of the Stanley-Reisner ideal
of $I_{A_n}$ then  there is an edge between $M$ and $N$ iff
$relint_{\bq}(M)\cap relint_{\bq}(N)\not= \{ 0 \}$. Every vector
in   $relint_{\bq}(M)$ can be written as a positive linear combination of the vectors 
$ {\bf e}_i$, where $i\in index(M)$. Since the vectors $\{{\bf e}_i|1\leq i\leq n\}$
are linearly independent, we conclude that  index($M)=$index($N$).\\
Therefore two minimal generators can be vertices of a connected
component of the graph $G_{pos_{\bq}(A_n)}$ if they have the same index.
The index of a minimal generator can be a set with three elements $\{i,j,k\}$ or
four elements $\{i,j,k,l\}$. By explicitly computing the edges among the 9 vertices  
with index $\{i,j,k\}$ we get that all of them are in the same connected component
which has 15 edges and looks like the FIGURE 1. Similarly, by explicitly computing 
the edges among the 12 vertices  
with index $\{i,j,k,l\}$ we get that all of them are in the same connected component
which has 18 edges and looks like the FIGURE 2.
\\
\bigskip
\newline
\setlength{\unitlength}{3947sp}%
\begingroup\makeatletter\ifx\SetFigFont\undefined%
\gdef\SetFigFont#1#2#3#4#5{%
  \reset@font\fontsize{#1}{#2pt}%
  \fontfamily{#3}\fontseries{#4}\fontshape{#5}%
  \selectfont}%
\fi\endgroup%
\begin{picture}(6716,2757)(850,-3400)
\thinlines
\put(7351,-1261){\circle*{100}}
\put(1501,-2901){\circle*{100}}
\put(2701,-2901){\circle*{100}}
\put(901,-1861){\circle*{100}}
\put(3301,-1861){\circle*{100}}
\put(1501,-821){\circle*{100}}
\put(1801,-2386){\circle*{100}}
\put(1801,-1336){\circle*{100}}
\put(2626,-1861){\circle*{100}}
\put(5101,-1861){\circle*{100}}
\put(5251,-1261){\circle*{100}}
\put(6301,-661){\circle*{100}}
\put(6901,-2911){\circle*{100}}
\put(6301,-3061){\circle*{100}}
\put(2701,-821){\circle*{100}}
\put(5751,-811){\circle*{100}}
\put(7351,-2461){\circle*{100}}
\put(5251,-2461){\circle*{100}}
\put(7501,-1861){\circle*{100}}
\put(6901,-811){\circle*{100}}
\put(5701,-2911){\circle*{100}}
\put(6901,-811){\makebox(1.6667,11.6667){\SetFigFont{5}{6}{\rmdefault}{\mddefault}{\updefault}.}}
\put(6901,-2911){\line( 1, 4){414.706}}
\put(1501,-811){\line( 3,-5){1244.118}}
\put(3301,-1861){\line(-1, 0){2400}}
\put(1801,-1336){\line( 0,-1){1050}}
\put(1801,-2386){\line( 3, 2){813.462}}
\put(2626,-1861){\line(-3, 2){813.462}}
\put(5101,-1861){\line( 1,-1){1200}}
\put(1501,-2911){\line( 3, 5){1239.706}}
\put(3301,-1861){\line(-3, 5){617.647}}
\put(2701,-821){\line(-1, 0){1200}}
\put(1501,-821){\line(-3,-5){617.647}}
\put(901,-1861){\line( 3,-5){617.647}}
\put(1501,-2901){\line( 1, 0){1200}}
\put(2701,-2901){\line( 3, 5){617.647}}
\put(5751,-811){\line( 1,-1){1625}}
\put(5251,-2461){\line( 4, 1){2258.823}}
\put(5251,-1261){\line( 4, 1){1658.823}}
\put(6301,-661){\line(-1,-4){564.706}}
\put(5101,-1861){\line( 1, 4){150}}
\put(5251,-1261){\line( 6, 5){516.393}}
\put(5751,-811){\line( 4, 1){552.941}}
\put(6301,-661){\line( 4,-1){600}}
\put(6901,-811){\line( 1,-1){450}}
\put(7351,-1261){\line( 1,-4){150}}
\put(7501,-1861){\line(-1,-4){150}}
\put(7351,-2461){\line(-1,-1){450}}
\put(6901,-2911){\line(-4,-1){600}}
\put(6301,-3061){\line(-4, 1){600}}
\put(5701,-2911){\line(-1, 1){450}}
\put(5251,-2461){\line(-1, 4){150}}
\put(3301,-1861){\line(-1, 0){2400}}
\put(7576,-1911){\makebox(0,0)[lb]{\smash{\SetFigFont{12}{14.4}{\rmdefault}{\mddefault}{\updefault}
(ik)(lj)
}}}
\put(5851,-3661){\makebox(0,0)[lb]{\smash{\SetFigFont{12}{14.4}{\rmdefault}{\mddefault}{\updefault}
FIGURE 2
}}}
\put(1801,-3661){\makebox(0,0)[lb]{\smash{\SetFigFont{12}{14.4}{\rmdefault}{\mddefault}{\updefault}
FIGURE 1
}}}
\put(2476,-3136){\makebox(0,0)[lb]{\smash{\SetFigFont{12}{14.4}{\rmdefault}{\mddefault}{\updefault}
(kj)(ik)
}}}
\put(1276,-3136){\makebox(0,0)[lb]{\smash{\SetFigFont{12}{14.4}{\rmdefault}{\mddefault}{\updefault}
(ki)(ij)
}}}
\put(1276,-736){\makebox(0,0)[lb]{\smash{\SetFigFont{12}{14.4}{\rmdefault}{\mddefault}{\updefault}
(ij)(jk)
}}}
\put(2476,-736){\makebox(0,0)[lb]{\smash{\SetFigFont{12}{14.4}{\rmdefault}{\mddefault}{\updefault}
(ji)(kj)
}}}
\put(3381,-1936){\makebox(0,0)[lb]{\smash{\SetFigFont{12}{14.4}{\rmdefault}{\mddefault}{\updefault}
(jk)(ki)
}}}
\put(270,-1936){\makebox(0,0)[lb]{\smash{\SetFigFont{12}{14.4}{\rmdefault}{\mddefault}{\updefault}
(ik)(ji)
}}}
\put(1831,-2521){\makebox(0,0)[lb]{\smash{\SetFigFont{12}{14.4}{\rmdefault}{\mddefault}{\updefault}
(jk)(ik)
}}}
\put(1831,-1261){\makebox(0,0)[lb]{\smash{\SetFigFont{12}{14.4}{\rmdefault}{\mddefault}{\updefault}
(ji)(ki)
}}}
\put(2551,-2061){\makebox(0,0)[lb]{\smash{\SetFigFont{12}{14.4}{\rmdefault}{\mddefault}{\updefault}
(kj)(ij)
}}}
\put(4320,-1911){\makebox(0,0)[lb]{\smash{\SetFigFont{12}{14.4}{\rmdefault}{\mddefault}{\updefault}
(il)(jk)
}}}
\put(6076,-3286){\makebox(0,0)[lb]{\smash{\SetFigFont{12}{14.4}{\rmdefault}{\mddefault}{\updefault}
(kl)(ji)
}}}
\put(6076,-586){\makebox(0,0)[lb]{\smash{\SetFigFont{12}{14.4}{\rmdefault}{\mddefault}{\updefault}
(ki)(lj)
}}}
\put(6976,-831){\makebox(0,0)[lb]{\smash{\SetFigFont{12}{14.4}{\rmdefault}{\mddefault}{\updefault}
(lk)(ji)
}}}
\put(5000,-831){\makebox(0,0)[lb]{\smash{\SetFigFont{12}{14.4}{\rmdefault}{\mddefault}{\updefault}
(ij)(kl)
}}}
\put(4500,-1281){\makebox(0,0)[lb]{\smash{\SetFigFont{12}{14.4}{\rmdefault}{\mddefault}{\updefault}
(ik)(jl)
}}}
\put(4500,-2536){\makebox(0,0)[lb]{\smash{\SetFigFont{12}{14.4}{\rmdefault}{\mddefault}{\updefault}
(ij)(lk)
}}}
\put(7426,-2536){\makebox(0,0)[lb]{\smash{\SetFigFont{12}{14.4}{\rmdefault}{\mddefault}{\updefault}
(il)(kj)
}}}
\put(7426,-1281){\makebox(0,0)[lb]{\smash{\SetFigFont{12}{14.4}{\rmdefault}{\mddefault}{\updefault}
(li)(jk)
}}}
\put(4900,-2986){\makebox(0,0)[lb]{\smash{\SetFigFont{12}{14.4}{\rmdefault}{\mddefault}{\updefault}
(kj)(li)
}}}
\put(7051,-2986){\makebox(0,0)[lb]{\smash{\SetFigFont{12}{14.4}{\rmdefault}{\mddefault}{\updefault}
(ki)(jl)
}}}
\put(1501,-2911){\line( 3, 5){1239.706}}
\put(7576,-1911){\makebox(0,0)[lb]{\smash{\SetFigFont{12}{14.4}{\rmdefault}{\mddefault}{\updefault}
(ik)(lj)
}}}
\end{picture}
\bigskip
\newline
Therefore we conclude that the graph $G_{pos_{\bq}(A_n)}$ has $\pmatrix{n \cr 3\cr}$ connected
components like the one in FIGURE 1, with  9 vertices and 15 edges each,  and
$\pmatrix{n \cr 4\cr}$ connected components like the one in FIGURE 2, with  12 vertices and 18 edges each.

\begin{Corollary} Let $L$ be a lattice with associated cone rationally affine equivalent to
$pos_{\bq}(A_n)$, then for the ideal $I_{L, \rho}$ we have that 
$$bar(I_{L, \rho})\ge 5\pmatrix{n \cr 3\cr}+6\pmatrix{n \cr 4\cr} \ and$$
$$ara_{A}(I_{L, \rho})\ge 4\pmatrix{n \cr 3\cr}+6\pmatrix{n \cr 4\cr}.$$
\end{Corollary}
\setlength{\unitlength}{3947sp}%
\begingroup\makeatletter\ifx\SetFigFont\undefined%
\gdef\SetFigFont#1#2#3#4#5{%
  \reset@font\fontsize{#1}{#2pt}%
  \fontfamily{#3}\fontseries{#4}\fontshape{#5}%
  \selectfont}%
\fi\endgroup%
\begin{picture}(6766,2757)(843,-3361)
\thinlines
\put(1801,-2386){\circle*{100}}
\put(7351,-2461){\circle*{100}}
\put(5251,-2461){\circle*{100}}
\put(1501,-821){\circle*{100}}
\put(901,-1861){\circle*{100}}
\put(2701,-821){\circle*{100}}
\put(3301,-1861){\circle*{100}}
\put(2701,-2901){\circle*{100}}
\put(1501,-2901){\circle*{100}}
\put(1801,-1336){\circle*{100}}
\put(2626,-1861){\circle*{100}}
\put(5251,-1261){\circle*{100}}
\put(6301,-661){\circle*{100}}
\put(5051,-1861){\circle*{100}}
\put(5751,-811){\circle*{100}}
\put(6901,-811){\circle*{100}}
\put(7351,-1261){\circle*{100}}
\put(6901,-2911){\circle*{100}}
\put(6301,-3061){\circle*{100}}
\put(7551,-1861){\circle*{100}}
\put(5701,-2911){\circle*{100}}
\put(6901,-811)
{\makebox(1.6667,11.6667){\SetFigFont{5}{6}{\rmdefault}{\mddefault}{\updefault}.}}
\put(5751,-811){\line( 1,-1){1625}}
\put(2701,-821){\line( 3,-5){617.647}}
\put(2701,-2901){\line(-1, 0){1200}}
\put(901,-1861){\line( 3, 5){617.647}}
\put(1801,-1336){\line( 3,-2){813.462}}
\put(2626,-1861){\line(-3,-2){813.462}}
\multiput(1801,-2386)(0.00000,123.52941){9}{\line( 0, 1){ 61.765}}
\put(5051,-1861){\line( 1, 3){200}}
\put(6901,-811){\line( 1,-1){450}}
\put(6901,-2911){\line(-4,-1){600}}
\put(5251,-2461){\line( 4, 1){2258.823}}
\put(6301,-661){\line(-1,-4){564.706}}
\put(5851,-3361){\makebox(0,0)[lb]
{\smash{\SetFigFont{12}{14.4}{\rmdefault}{\mddefault}{\updefault}
FIGURE 4
}}}
\put(1721,-3361){\makebox(0,0)[lb]{\smash{\SetFigFont{12}{14.4}
{\rmdefault}{\mddefault}{\updefault}
FIGURE 3
}}}
\end{picture}

\demo  Recall that $b_G$ is the smallest number $s$ of
 subgraphs $B_i$ of $G$,
consisting of two vertices and an edge or just a vertex, such that the graph
$\cup_{i=1}^sB_i$ is spanning. For the $\pmatrix{n \cr 3\cr}$ 
connected components of $G_{pos_{\bq}(A_n)}$,
 like the one in
FIGURE 1, this number is five as it can be seen in FIGURE 3. While for the 
$\pmatrix{n \cr 4\cr}$ connected 
components of $G_{pos_{\bq}(A_n)}$, like the one in
FIGURE 2, this number is six as it can be seen in FIGURE 4. Thus $b_{G_{pos_{\bq}(A_n)}}=
5\pmatrix{n \cr 3\cr}+6\pmatrix{n \cr 4\cr}$. \\
Recall also that $c_G$ is the smallest number $s$ of 
complete subgraphs $G_i$ of $G$, such that the graph
$\cup_{i=1}^sG_i$ is spanning. Note also that graphs like those in FIGURE 1 have only
one complete subgraph  with 3 vertices 
and  those in FIGURE 2 have only complete subgraphs
with two or one vertices. Consequently, for the 
$\pmatrix{n \cr 3\cr}$ connected components
of $G_{pos_{\bq}(A_n)}$, like the one in
FIGURE 1, the number $c_G$ is four as it can be seen in FIGURE 3. 
While for the $\pmatrix{n \cr 4\cr}$ connected 
components,  like the one in
FIGURE 2, this number is six as it can be seen in FIGURE 4. 
Therefore $c_{G_{pos_{\bq}(A_n)}}=
4\pmatrix{n \cr 3\cr}+6\pmatrix{n \cr 4\cr}$. \\
The proof follows from Theorem 4.5.
\bigskip
\newline
\par Next we will prove that the lower bounds computed in Corollary 5.2 are sharp by computing the exact value of the binomial arithmetical rank and the A-homogeneous arithmetical rank for the toric ideal $I_{A_n}$.
\begin{proposition} The ideal $I_{A_n}$ is generated up to radical  
by the $5\pmatrix{n \cr 3\cr}+6\pmatrix{n \cr 4\cr}$ binomials 
$ x_{ij}x_{kj}-x_{jk}x_{ik}, 
x_{ij}^2x_{ki}-x_{ik}^2x_{ji}, x_{ij}x_{kl}-x_{il}x_{kj}$, where $i,j,k,l \in \{1,\ldots,n\}$. Therefore 
$bar(I_{A_n})= 5\pmatrix{n \cr 3\cr}+6\pmatrix{n \cr 4\cr}.$
\end{proposition}
\demo  Let $J$ be the ideal in $K[\{x_{ij}|1\leq i,j\leq n,i\not= j\}]$ 
generated by
the $5\pmatrix{n \cr 3\cr}+6\pmatrix{n \cr 4\cr}$ binomials 
$x_{ij}x_{kl}-x_{il}x_{kj}, x_{ij}x_{kj}-x_{jk}x_{ik}, x_{ij}^2x_{ki}-x_{ik}^2x_{ji}$, 
where $i,j,k,l \in \{1,\ldots,n\}$. We will use Hilbert's Nullstellensatz 
to prove the theorem. Obviously $J\subset I_{A_n}$ and 
therefore $V(I_{A_n})\subset V(J)$. Note that the toric variety $V(I_{A_n})$
is the Zariski-closure of the point $P=(1,\dots ,1)\in {K}^{n(n-1)}$ under the toric
action induced by the set of vectors $A_n$.\\
Let ${\bf a}\in {K}^{n(n-1)}$ be a point in $V(J)$ with $a_{ij}\not= 0$, 
for some fixed indices $i,j$. There are two cases:\\
a) $a_{ji}=0$. Then, using the binomials $x_{ij}^2x_{ki}-x_{ik}^2x_{ji}$ 
and $x_{ji}x_{ki}-x_{ij}x_{kj}$, we get
 that $a_{ki}=0$ and $a_{kj}=0$ for every index $k$ different from $i,j$. \\
In addition, using the binomials $x_{ij}x_{kl}-x_{il}x_{kj}$ and $x_{ji}^2x_{kj}-x_{jk}^2x_{ij}$,
we have that $a_{kl}=0$ and $a_{jk}=0$ for every indices $k,l$ different from $i,j$. \\
Let $T=\{k|a_{ik}\not= 0\}$. Note that $T$ is not empty, because $j\in T$. 
Let $F_{i,T}=pos_{\bq}\{2{\bf e}_i+{\bf e}_k/k\in T\}$, then $F_{i,T}$ is a face 
of $\sigma $ whose defining vector is ${\bf c}_{i,T}$.
Setting $t_i=1$, $t_k=a_{ik}$, for every $k\in T$, and $t_l=0$, for every 
$l\notin T$, we obtain that ${\bf a}$ is in the orbit of the point $P_{F_{i,T}}$. 
Thus ${\bf a}$ belongs to $V(I_{A_n})$.\\
b)$a_{ji}\not=0$. Let $T=\{k|a_{ik}\not=0\}\cup \{i\}$. Note that $j\in T$. Let $k\in T$ 
then, from the definition, $a_{ik}\not=0$. Using the binomial 
$x_{ij}^2x_{ki}-x_{ik}^2x_{ji}$ we obtain that $a_{ki}\not=0$. 
Then, from the binomial $x_{kj}x_{ij}-x_{ji}x_{ki}$ we have 
that $a_{kj}\not=0$. Finally, from the binomial $x_{ij}x_{kj}-x_{jk}x_{ik}$,
 we
conclude that $a_{jk}\not=0$.\\
Let $\{k,l\}\subset T$ and $\{k,l\}\cap \{i,j\}=\emptyset $, then, using the binomial $x_{ij}x_{kl}-
x_{il}x_{kj}$, we take $a_{kl}\not=0$. \\
Assume that $k\notin T$, 
then, from the definition, $a_{ik}=0$. The binomial 
$x_{ij}^2x_{ki}-x_{ik}^2x_{ji}$ gives $a_{ki}=0$, while the 
binomial $x_{kj}x_{ij}-x_{ji}x_{ki}$ gives $a_{kj}=0$. From the binomial $x_{ji}^2x_{kj}-x_{jk}^2x_{ij}$ we
conclude that $a_{jk}=0$. Also $a_{kl}=0$ for every index $l$, 
because of the binomial $x_{ij}x_{kl}-x_{il}x_{kj}$. 
The binomial $x_{ij}x_{lk}-x_{ik}x_{lj}$ give us that $a_{lk}=0$ for every index $l$.\\
Therefore $a_{pq}\not=0$ if $\{p,q\}\subset T$, while $a_{pq}=0$ if $\{p,q\}\not \subset T$.
Let $F_T=pos_{\bq}(2{\bf e}_p+{\bf e}_q/ \{p,q\}\subset T)$, then $F_T$ is a face of $\sigma $ 
whose defining vector is ${\bf c}_T$.\\ 
We will prove that the point ${\bf a}$ is in the orbit of the point $P_{{F}_T}$.
Let $i,j\in T$ and $\omega $ be any cubic root of $a_{ij}a_{ji}$. 
Setting $t_i=a_{ij}\omega ^{-1}$
and $t_j=a_{ji}\omega ^{-1}$ we have $a_{ij}=t_i^2t_j$ and $a_{ji}=t_j^2t_i$. 
For any $k\in T$ put $t_k=a_{jk}t_j^{-2}$, then of course $a_{jk}=t_j^2t_k$.\\
Using the binomials $x_{ji}^2x_{kj}-x_{jk}^2x_{ij}$, $x_{ij}x_{kj}-x_{jk}x_{ik}$, 
$x_{ij}^2x_{ki}-x_{ik}^2x_{ji}$ we conclude step by step that $a_{kj}=t_k^2t_j$,
$a_{ik}=t_i^2t_k$ and $a_{ki}=t_k^2t_i$. Then for any two $k,l$ in $T$, from the binomial
 $x_{ij}x_{kl}-x_{il}x_{kj}$, we have that $a_{kl}=t_k^2t_l$. 
Put $t_l=0$ for all $l\notin T$. Then the point ${\bf a}$ is in the orbit of the point $P_{{F}_T}$, so it is a point of $V(I_{A_n})$.\\
The second part of the proposition now follows from corollary 5.2.

\begin {remark} {\rm We can choose the binomials $x_{ij}^2x_{jk}-x_{ji}^2x_{ik}$ 
instead of $x_{ij}^2x_{ki}-x_{ik}^2x_{ji}$ to generate the radical of $I_{A_n}$. 
In addition, from the proof
of the above theorem we can see that the faces of the cone $\sigma $ are in the form
$F_{i,T}$ or $F_T$, for all the possible choises of $i$ and $T$.} 
\end {remark} 

\medskip
\begin{proposition} The $A_n$-homogeneous arithmetical rank of $I_{A_n}$ is equal to 
$4\pmatrix{n \cr 3\cr}+6\pmatrix{n \cr 4\cr}$. 
\end{proposition}
\demo  The ideal $I_{A_n}$ is generated up to radical 
by the $A_n$-homogeneous polynomials
$x_{ij}^2x_{ki}-x_{ik}^2x_{ji}, 
x_{ij}^3x_{kj}^3-x_{jk}^3x_{ik}^3+x_{ki}^3x_{ji}^3-x_{ij}^2x_{ki}^2x_{jk}^2, 
 x_{ij}x_{kl}-x_{il}x_{kj} $, where $i,j,k,l \in \{1,\ldots ,n\}$. The proof
follows from proposition 5.3 and the observation that $ (x_{ij}x_{kj}-x_{jk}x_{ik})^5$
belongs to the ideal generated by the previous $A_n$-homogeneous polynomials. 
Note that the $3\pmatrix{n \cr 3\cr}$ binomials $x_{ij}^2x_{ki}-x_{ik}^2x_{ji}$
correspond to the $3\pmatrix{n \cr 3\cr}$ complete subgraphs of $G_{pos_{\bq}(A_n)}$
with two vertices
like those in FIGURE 3. The $\pmatrix{n \cr 3\cr}$ $A_n$-homogeneous polynomials 
$x_{ij}^3x_{kj}^3-x_{jk}^3x_{ik}^3+x_{ki}^3x_{ji}^3-x_{ij}^2x_{ki}^2x_{jk}^2$ 
correspond to the $\pmatrix{n \cr 3\cr}$ 
complete subgraphs of $G_{pos_{\bq}(A_n)}$ with three vertices
like the one in FIGURE 3. The $6\pmatrix{n \cr 4\cr}$ binomials $x_{ij}x_{kl}-x_{il}x_{kj} $
correspond to the $6\pmatrix{n \cr 4\cr}$ complete subgraphs of $G_{pos_{\bq}(A_n)}$
with two vertices
like those in FIGURE 4.

\begin {remark} {\rm The bounds given in corollary 5.2 are also bounds
for the minimal number of generators of a lattice ideal $I_{L, \rho}$
with associated cone rationally affine equivalent to
$pos_{\bq}(A_n)$. In particular for any such ideal 
the minimal generators are greater 
than or equal to $5\pmatrix{n \cr 3\cr}+6\pmatrix{n \cr 4\cr}$. 
This implies that for any such ideal, for $n\ge 3$,
it is impossible  to be complete intersection, since $h(I_{L, \rho})=n(n-2)$.} 
\end {remark} 
\begin {remark} {\rm While theorem 4.1 give lower bounds
for $ara_A(I_{L,\rho })$, $bar(I_{L,\rho })$ and $\mu(I_{L,\rho })$ it does not provide 
a lower bound for $ara(I_{L,\rho })$. Nevertheless thorem 4.1 provides certain information
on the form and size of the polynomials $F_1,\ldots ,F_s$ such that
 $rad(I_{L,\rho })=rad(F_1,\ldots ,F_s)$. We know that for every vertex we need atleast
one monomial in atleast one of the $F_1,\ldots ,F_s,$  corresponding to the vertex. 
In particular for the ideals  $I_{L,\rho }$ studied in this section
we know that $n(n-2)\leq ara(I_{L,\rho })\leq n(n-1)$, by the Krull's principal ideal
theorem and the results of Eisenbud, Evans and Storch \cite{E-E}, \cite{Sto}. From 
theorem 4.1 we know that in these $s=ara(I_{L,\rho })$ polynomials there must be
atleast $9\pmatrix{n \cr 3\cr}+12\pmatrix{n \cr 4\cr}$ monomials, in atleast
$4\pmatrix{n \cr 3\cr}+6\pmatrix{n \cr 4\cr}$ $A_n$-homogeneous components. 
For example for $n=10$ we know that we need a number of polynomials between $80$ to $90$
to generate the radical of, say, $I_{A_{10}}$. Those polynomials should have totally
atleast 3600 monomials, so on the average atleast $40$ to $45$, and 
therefore even for small $n$'s 
the polynomials involved are huge. It will be an interesting problem to compute these
polynomials even for $n=10$.} \end {remark} 

\par Note that in all the cases, that we know explicitly the polynomials which define
a lattice ideal up to radical, the polynomials involved are all $A$-homogeneous, 
see \cite{B-M-T1}, \cite{R-V}, \cite{Th}. The results of this paper show
that $A$-homogeneous polynomials are not always enough to define 
a lattice ideal up to radical.
Therefore we have to understand better the topic of non $A$-homogeneous 
set theoretic intersections for lattice ideals.  
Also these results give a different perspective relative to the 
famous Macaulay curve $(t^4,t^3u,tu^3,u^4)$ 
in the three
dimensional projective space, for which we know that it is not 
A-homogeneous set-theoretic complete intersection, see \cite{Tho}.

\par {\sc Department of Mathematics,
University of Ioannina, Ioannina 45110 (GREECE)}\\ 

{\sc Universit\'e de Grenoble I, Institut Fourier,
UMR 5582, B.P.74, 38402 Saint-Martin D'H\`eres Cedex, and IUFM de Lyon, 5 rue Anselme, 69317 Lyon Cedex (FRANCE)}\\ 

{\sc Department of Mathematics,
University of Ioannina, Ioannina 45110 (GREECE)}
\end{document}